\documentclass[12pt]{article}
\usepackage{amsfonts}
\usepackage{mathrsfs}

\usepackage{amsmath,amssymb,amsfonts,theorem,makeidx,latexsym,epsfig,color,subfigure,algorithm,algorithmic}\usepackage[numbers,sort&compress]{natbib}

\newtheorem{defn}{Definition}[section]

\newtheorem{lem}[defn]{Lemma}

{\theorembodyfont{\rmfamily}

}

\newtheorem{thm}[defn]{Theorem}

\newtheorem{prop}[defn]{Proposition}

\newtheorem{obs}[defn]{Observation}

\numberwithin{equation}{section}

\def\bp{{\noindent\bf Proof. \ }}

\def\ep{\hfill$\square$\par\bigskip}

\title{Super domination in trees \footnote{The research is supported by NSFC (No. 11301440),
Natural Science Foundation of Fujian Province (CN)(2015J05017)}}

\author{{\large Zhuang Wei\thanks{Corresponding author; E-mail: zhuangweixmu@163.com}} \\
{\it \normalsize School of Applied Mathematics},
\\{\it \normalsize  Xiamen University of Technology, Xiamen 361024, P.R.China }}

\date{}

\begin{document}

\maketitle

\begin{abstract} For $S\subseteq V(G)$, we define $\overline{S}=V(G)\setminus
S$. A set $S\subseteq V(G)$ is called a super dominating set if for
every vertex $u\in \overline{S}$, there exists $v\in S$ such that
$N(v)\cap \overline{S}=\{u\}$. The super domination number
$\gamma_{sp}(G)$ of $G$ is the minimum cardinality among all super
dominating sets in $G$. The super domination subdivision number
$sd_{\gamma_{sp}}(G)$ of a graph $G$ is the minimum number of edges
that must be subdivided in order to increase the super domination
number of $G$. In this paper, we investigate the ratios between
super domination and other domination parameters in trees. In
addition, we show that for any nontrivial tree $T$, $1\leq
sd_{\gamma_{sp}}(T)\leq 2$, and give constructive characterizations
of trees whose super domination subdivision number are $1$ and $2$,
respectively.
\end{abstract}

\begin{minipage}{150mm}

\end{minipage}

\section{Introduction}
Let $G=(V, E)$ be a simple graph without isolated vertices, and let
$v$ be a vertex in $G$. The \emph{open neighborhood} of $v$ is
$N(v)=\{u\in V|uv\in E\}$ and the \emph{closed neighborhood} of $v$
is $N[v]=N(v)\cup \{v\}$. For $S\subseteq V(G)$, we define
$\overline{S}=V(G)\setminus S$. The \emph{degree} of a vertex $v$ is
$d(v)=|N(v)|$. For two vertices $u$ and $v$ in a connected graph
$G$, the \emph{distance} $d(u, v)$ between $u$ and $v$ is the length
of a shortest $(u, v)$-path in $G$. The maximum distance among all
pairs of vertices of $G$ is the \emph{diameter} of a graph $G$ which
is denoted by $diam(G)$. A \emph{leaf} of $G$ is a vertex of degree
$1$, and a \emph{support vertex} of $G$ is a vertex adjacent to a
leaf. A support vertex that is adjacent to at least two leaves we
call a \emph{strong support vertex}. The corona $G\circ K_1$ is the
graph obtained from a graph $G$ by attaching a leaf to each vertex
$v\in V(G)$.

A \emph{dominating set} (respectively, \emph{total dominating set})
in a graph $G$ is a set $S$ of vertices of $G$ such that every
vertex in $V(G)\setminus S$ (respectively, $V(G)$) is adjacent to at
least one vertex in $S$. The \emph{domination number} (respectively,
\emph{total domination number}) of $G$, denoted by $\gamma(G)$
(respectively, $\gamma_t(G)$), is the minimum cardinality of a
dominating set (respectively, total dominating set) of $G$. A
dominating set (respectively, total dominating set) of $G$ with
cardinality $\gamma(G)$ (respectively, $\gamma_t(G)$) is called a
\emph{$\gamma(G)$-set} (respectively, \emph{$\gamma_t(G)$-set}). We
say a vertex $v$ in $G$ is total dominated, by a set $D$, if
$N(v)\cap D\neq \emptyset$.

The study of super domination in graphs was introduced in
\cite{Lemanska}. A set $S\subseteq V(G)$ is called a \emph{super
dominating set} if for every vertex $u\in \overline{S}$, there
exists $v\in S$ such that $N(v)\cap \overline{S}=\{u\}$. In
particular, we say that $v$ is an external private neighbor of $u$
with respect to $\overline{S}$. For a super dominating set $S$ of
$G$, let $P_S(G)=\{v|$ $v$ is an external private neighbor of $u$
with respect to $\overline{S}$, for each $u\in \overline{S}\}$,
$Q_S(G)=\{v|v$ belongs to $\overline{S}$ and $v$ has only one
external private neighbor with respect to $\overline{S}\}$ and
$U_S(G)=\{v|$ $v$ is the unique external private neighbor of $u$
with respect to $\overline{S}$, for each $u\in Q_S(G)\}$. The
\emph{super domination number} $\gamma_{sp}(G)$ of $G$ is the
minimum cardinality among all super dominating sets in $G$. A super
dominating set of $G$ with cardinality $\gamma_{sp}(G)$ is called a
\emph{$\gamma_{sp}(G)$-set}. More results in this area were
investigated in \cite{Dettlaff, Krishnakumari} and elsewhere.

The domination subdivision number $sd_\gamma(G)$ of a graph $G$ is
the minimum number of edges that must be subdivided (where each edge
in $G$ can be subdivided at most once) in order to increase the
domination number. The domination subdivision number was first
introduced in Velammal's thesis \cite{Velammal} and since then many
results have also been obtained on the parameters $sd_\gamma(G)$,
$sd_{\gamma_t}(G)$, $sd_{\gamma_2}(G)$, $sd_{\gamma_{pr}}(G)$ (see
\cite{Atapour, Favaron, Haynes3, Karami, Favaron2, Atapour2}). One
of the purpose of this paper is to initialize the study of the super
domination subdivision number. The super domination subdivision
number $sd_{\gamma_{sp}}(G)$ of a graph $G$ is the minimum number of
edges that must be subdivided in order to increase the super
domination number of $G$ (each edge in $G$ can be subdivided at most
once).

In this paper, we investigate the ratios between super domination
and other domination parameters in trees. In addition, we show that
for any nontrivial tree $T$, $1\leq sd_{\gamma_{sp}}(T)\leq 2$, and
give constructive characterizations of trees whose super domination
subdivision number are $1$ and $2$, respectively.

\section{On the ratios between super domination and other domination parameters in trees}

From the definitions of domination number, total domination number
and super domination number, we have the following observations.

\begin{obs}
Let $G$ be a connected graph that is not a star. Then,

$(1)$ there is a $\gamma$-set of $G$ that contains no leaf, and

$(2)$\cite{Henning1} there is a $\gamma_{t}$-set of $G$ that
contains no leaf.
\end{obs}

\begin{obs}
Let $G$ be a connected graph of order at least $2$, $v$ be a support
vertex of $G$ and $S$ be a $\gamma_{sp}$-set of $G$. Then, at most
one of $v$ and its leaf-neighbors belongs to $\overline{S}$.
\end{obs}

\begin{obs}
Let $T$ be a tree containing the strong support vertices $u_1, u_2,
\cdots, u_t$, and $T'$ be a tree is obtained from $T$ by deleting
$x_i$ $(0\leq x_i \leq d(u_i)-2)$ leaf-neighbors of each $u_i$,
$i=1, 2, \cdots, t$. Then, $\gamma_{sp}(T)=a$ if and only if
$\gamma_{sp}(T')=a-\sum \limits_{i=1}^{t} x_i$.
\end{obs}

\begin{prop}
Let $T$ be a tree of order at least $2$ and $v$ be a leaf of $T$,
there is always a $\gamma_{sp}$-set $S$ of $T$ such that $v\in
\overline{S}$.
\end{prop}

Next, we will investigate how to obtain a $\gamma_{sp}$-set of $T$
mentioned in Proposition~2.4. Given an arbitrary
$\gamma_{sp}(T)$-set $S$, we root the tree $T$ at the leaf $v$. For
any vertex $u$, let $C(u)$ be the set consisting of the children of
$u$. Distinguish three cases as follows:

\vspace{0.3cm}

$(I)$ $v\in \overline{S}$.

In this case, the set $S$ is the desired $\gamma_{sp}$-set of $T$.

\vspace{0.3cm}

$(II)$ $v\not \in \overline{S}$ and its support vertex, say $u$,
belongs to $\overline{S}$.

In this case, $v$ is the external private neighbor of $u$ with
respect to $\overline{S}$. We will determine a set $H$ by the
following procedure:

\vspace{0.3cm}

$(1)$ We set $H:=(\overline{S}\setminus \{u\})\cup \{v\}$,
$X_0:=C(u)\cap \overline{S}$ and $i:=0$. Take a vertex $t\in X_i$,
and set $P:=\{t\}$.

\vspace{0.2cm}

$(2)$ We query whether $X_i$ is an empty set or not.

--- If the answer to the query is `yes' and $i=0$,

then we terminate.

--- If the answer to the query is `yes' and $i\neq 0$,

then we set $i:=i-1$.

Set $P=\emptyset$, take a vertex $t\in X_i$ and put it into $P$. Go
to $(2)$.

--- If the answer to the query is `no' and $P\cap \overline{S}\neq \emptyset$,

then go to $(3)$.

--- If the answer to the query is `no' and $P\cap \overline{S}=\emptyset$,

then go to $(4)$.

\vspace{0.2cm}

$(3)$ We take a vertex $x\in X_i$ and one of the external private
neighbors of $x$ with respect to $\overline{S}$, say $y$ (Note that
$y\in C(x)$).

Set $X_i:=X_i\setminus \{x\}$, $H:=(H\setminus \{x\})\cup \{y\}$.

Set $i:=i+1$, $X_i=C(x)\cap \overline{S}$.

Set $i:=i+1$, $X_i:=C(y)\cap P_S(T)$.

Set $P=\emptyset$, take a vertex $t\in X_i$ and put it into $P$. Go
to $(2)$.

 \vspace{0.2cm}

$(4)$ We take a vertex $w\in X_i$ and a vertex $z\in C(w)\cap
\overline{S}$ (Note that $w$ is an external private neighbors of $z$
with respect to $\overline{S}$).

Set $X_i:=X_i\setminus \{w\}$, $H:=(H\setminus \{z\})\cup \{w\}$.

Set $i:=i+1$, $X_i=C(w)\cap P_S(T)$.

Set $i:=i+1$, $X_i:=C(z)\cap \overline{S}$.

Set $P=\emptyset$, take a vertex $t\in X_i$ and put it into $P$. Go
to $(2)$.

\vspace{0.3cm}

After the end of this procedure, the set $V(T)\setminus H$ is the
desired $\gamma_{sp}$-set of $T$.

\vspace{0.3cm}

$(III)$ Both of $v$ and its support vertex do not belong to
$\overline{S}$.

In this case, the support vertex of $v$, say $u$, has at least one
child which belongs to $\overline{S}$.

If there exists a vertex $u_1\in C(u)\cap \overline{S}$ such that
$u$ is an external private neighbors of $u_1$ with respect to
$\overline{S}$, then let $D=(\overline{S}\setminus \{u_1\})\cup
\{v\}$, and the set $V(T)\setminus D$ is the desired
$\gamma_{sp}$-set of $T$.

If $u\not \in P_S(T)$, then it has at least two children belonging
to $\overline{S}$. Assume that $C(u)\cap P_S(T)=\emptyset$, then let
$D=\overline{S}\cup \{u\}$, and $V(T)\setminus D$ is a super
dominating set of $T$ whose cardinality is less than $|S|$, a
contradiction. Hence, we consider the case of $C(u)\cap P_S(T)\neq
\emptyset$. Let $C(u)\cap P_S(T)=\{u_1, u_2, \cdots, u_k\}$. For any
$i\in \{1, 2, \cdots, k\}$, let $T_i$ be the component of $T-uu_i$
containing $u_i$. Clearly, $S_i=S\cap V(T_i)$ is a
$\gamma_{sp}(T_i)$-set. Similar to the argument of Case $(II)$,
there must be a $\gamma_{sp}(T_i)$-set $S'_i$ such that $u_i\in
\overline{S'_i}$. Then, $(S\setminus (\bigcup \limits_{i=1}^{k}
S_i))\cup (\bigcup \limits_{i=1}^{k} S'_i)$ is also a
$\gamma_{sp}(T)$-set, and similar to the case of $C(u)\cap
P_S(T)=\emptyset$, leading to a contradiction.

\vspace{0.8cm}

\begin{center}
  \includegraphics[width=2.3in]{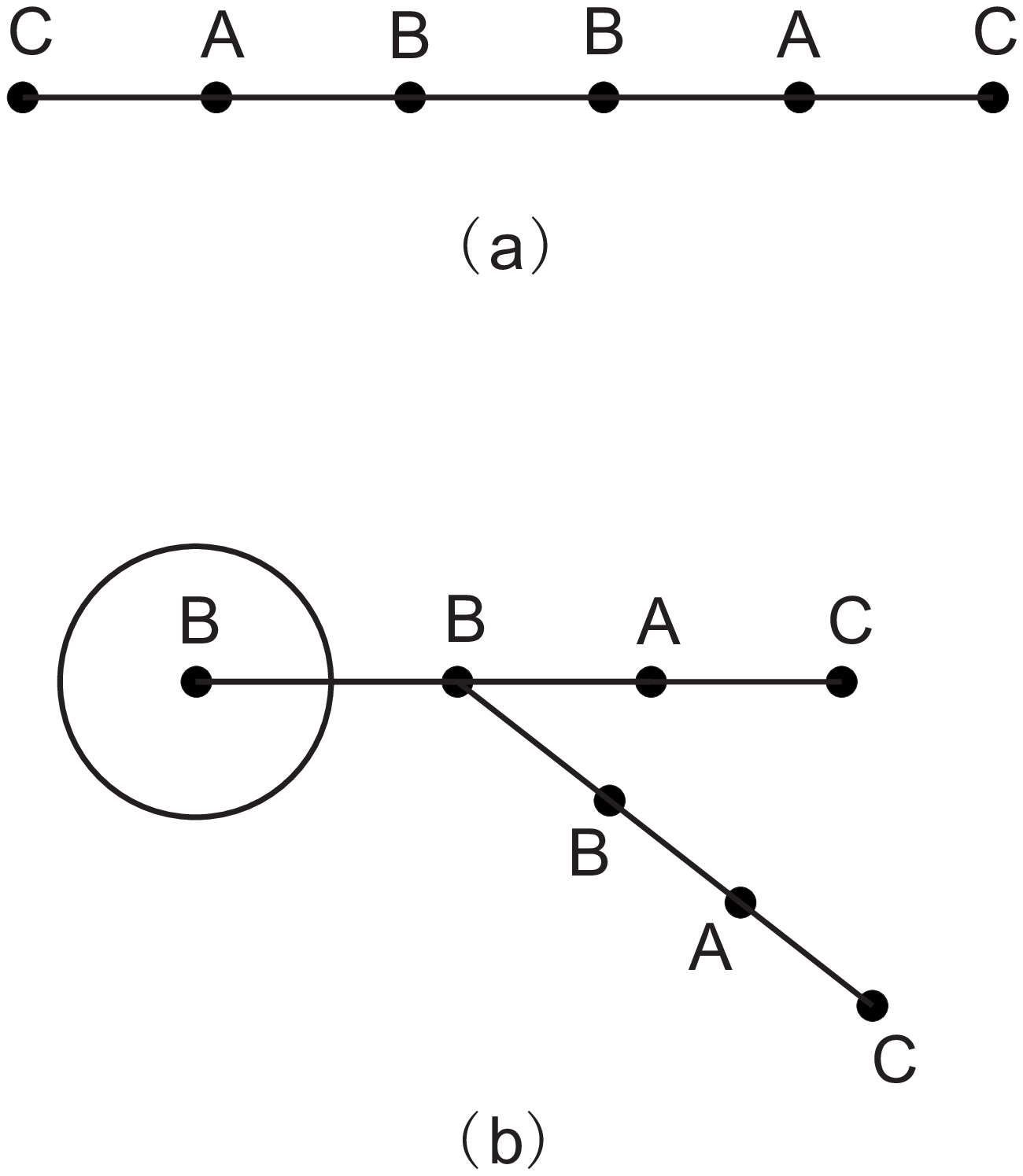}
 \end{center}
\qquad \qquad \qquad \qquad \qquad \qquad \qquad \qquad \qquad
{\small Fig.1} $\\$

By a weak partition of a set we mean a partition of the set in which
some of the subsets may be empty. We define a \emph{labeling} of a
tree $T$ as a weak partition $S=(S_A, S_B, S_C)$ of $V(T)$ (This
idea of labeling the vertices is introduced in \cite{Dorfling}). We
will refer to the pair $(T, S)$ as a \emph{labeled tree}. The label
or \emph{status} of a vertex $v$, denoted sta$(v)$, is the letter
$x\in \{A, B, C\}$ such that $v\in S_x$. Next, let $\mathscr{T}$ be
the family of labeled trees that: (i) contains $(P_6, S_0)$ where
$S_0$ is the labeling that assigns to the two leaves of the path
$P_6$ status $C$, to the support vertices status $A$ and to the
remaining vertices status $B$ (see Fig.1(a)); and (ii) is closed
under the operation $\mathscr{O}$ that is listed below, which extend
the tree $T'$ to a tree $T$ by attaching a $P_6$ to the vertex $v\in
V(T')$.

{\bf Operation} $\mathscr{O}$: Let $v$ be a vertex with sta$(v)=B$.
Add a path $u_1u_2u_3u_4u_5u_6$ and the edge $u_3v$. Let
sta$(u_1)=$sta$(u_6)=C$, sta$(u_2)=$sta$(u_5)=A$ and
sta$(u_3)=$sta$(u_4)=B$.

The operation $\mathscr{O}$ is illustrated in Fig.2(b).

For a tree $T$, we are ready to investigate the ratios between the
super domination number and two most important domination
parameters, namely the domination number and the total domination
number. We obtain the following conclusions.

\begin{thm}
For any nontrivial tree $T$, we have that
$0<\frac{\gamma(T)}{\gamma_{sp}(T)}\leq 1$. Further, the nontrivial
trees $T$ satisfying $\frac{\gamma(T)}{\gamma_{sp}(T)}=1$ are
precisely the corona $H\circ K_1$ for any tree $H$.
\end{thm}

\begin{thm}
Let $T$ be a tree of order at least $3$, we have that
$0<\frac{\gamma_t(T)}{\gamma_{sp}(T)}\leq \frac{4}{3}$. Further, the
trees $T$ of order at least $3$ satisfying
$\frac{\gamma_t(T)}{\gamma_{sp}(T)}=\frac{4}{3}$ are precisely those
trees $T$ such that $(T, S)\in \mathscr{T}$ for some labeling $S$.
\end{thm}

First, we show the lower bounds of Theorem~2.5 and Theorem~2.6 are
optimal.

Given any tree $T$, and construct a sequence of trees $T_0, T_1,
T_2, \cdots$, such that $T_0=T$, the tree $T_{i+1}$ is obtained from
$T_i$ by adding a vertex and joining it to one of the support
vertices of $T_i$, $i=0, 1, 2, \cdots$. By Observation~2.1 and 2.3,
both of the domination number and the total domination number of the
resulting trees have never changed in this process, but the super
domination number is constantly increasing. It means that when the
number $n$ is sufficiently large, the two ratios
$\frac{\gamma(T_n)}{\gamma_{sp}(T_n)}$ and
$\frac{\gamma_t(T_n)}{\gamma_{sp}(T_n)}$ are close to $0$.

It is well known that for any graph $G$ without isolated vertices,
we have that $\gamma(G)\leq \lfloor \frac{n}{2}\rfloor$ and
$\gamma_{sp}(G)\geq \lceil \frac{n}{2}\rceil$. It means that for any
tree $T$, $\gamma(T)\leq \gamma_{sp}(T)$.  The two theorems as
follows are useful to prove Theorem~2.5.

\begin{thm}\cite{Payan}
For a graph $G$ with even order $n$ and no isolated vertices,
$\gamma(G)=\frac{n}{2}$ if and only if the components of $G$ are the
cycle $C_4$ or the corona $H\circ K_1$ for any connected graph $H$.
\end{thm}

\begin{thm}\cite{Lemanska}
Let $T$ be a tree of order at least two. Then,
$\gamma_{sp}(T)=\frac{n}{2}$ if and only if $T\in \mathscr{R}$.
\end{thm}

The family $\mathscr{R}$ mentioned in Theorem~2.8 is a family of
trees that can be obtained from a sequence of trees $T_1, T_2,
\cdots, T_j (j\geq 1)$ such that:

$(1)$ The tree $T_1=P_2=a_1b_1$.

$(2)$ If $j\geq 2$, define the tree $T_j$ such that
$V(T_j)=V(T_{j-1})\cup \{a_j, b_j\}$ and $E(T_j)=E(T_{j-1})\cup
\{a_jb_j\}\cup \{e\}$, where $e=a_ia_j$ or $e=b_ib_j$ for some
$1\leq i\leq j-1$.

If $\frac{\gamma(T)}{\gamma_{sp}(T)}=1$, then
$\gamma(T)=\gamma_{sp}(T)=\frac{n}{2}$. It follows from Theorem~2.7
and Theorem~2.8 that $T$ is a corona $H\circ K_1$ for some tree $H$.
On the other hand, it is easy to see that for any tree $H$, both the
domination number and the super domination number of the corona
$H\circ K_1$ are $\frac{n}{2}$.

Finally, we are ready to prove the second half of Theorem~2.6. Let
$T$ be a tree containing strong support vertices, and $T'$ be the
tree obtained from $T$ by deleting all except one leaf-neighbor of
every strong support vertex of $T$. Then, it follows from
Observation~2.1 and 2.3 that $\gamma_t(T')=\gamma_t(T)$ and
$\gamma_{sp}(T')<\gamma_{sp}(T)$. And then
$\frac{\gamma_t(T)}{\gamma_{sp}(T)}<\frac{\gamma_t(T')}{\gamma_{sp}(T')}$.
So when we prove Theorem~2.6, we only need to consider the trees
containing no strong support vertex.

Before this, we present a preliminary result.

\begin{obs}
Let $T$ be a tree and $S$ be a labeling of $T$ such that $(T, S)\in
\mathscr{T}$. Then, $T$ has the following properties:

$(a)$ A vertex is labeled $C$ if and only if it is a leaf.

$(b)$ A vertex is labeled $A$ if and only if it is a support vertex.

$(c)$ Every support vertex has degree two, and its two neighbor have
status $C$ and $B$, respectively.

$(d)$ If a vertex has status $B$, then all of its neighbors have
status $B$ except one which has status $A$.

$(e)$ $|S_A|=|S_B|=|S_C|$.

$(f)$ $S_A\cup S_B$ is a $\gamma_t$-set of $T$.
\end{obs}

\begin{lem}
Let $T$ be a tree and $S$ be a labeling of $T$ such that $(T, S)\in
\mathscr{T}$. Then for any leaf $v$ of $T$, there exists a set $X$
of cardinality $\gamma_t(T)-1$ such that $v$ belongs to $X$, and
each vertex of $T$ is total dominated by $X$ except for $v$.
\end{lem}

\bp Let $v$ be a leaf of $T$, $u$ be its support vertex and $w$ be
the neighbor of $u$. It follows from Observation~2.9(b), (c) and (f)
that $u$ , $w$ have status $A$ and $B$, respectively, and $D=S_A\cup
S_B$ is a $\gamma_t$-set of $T$. It is easy to see that $(D\setminus
\{u, w\})\cup \{v\}$ is the set as we desired. \ep

\begin{lem}
Let $T$ be a tree and $S$ be a labeling of $T$ such that $(T, S)\in
\mathscr{T}$. Then, $\frac{\gamma_t(T)}{\gamma_{sp}(T)}=
\frac{4}{3}$.
\end{lem}

\bp The proof is by induction on the number $h(T)$ of operations
required to construct the tree $T$. Observe that $T=P_6$ when
$h(T)=0$, and clearly $\frac{\gamma_t(T)}{\gamma_{sp}(T)}=
\frac{4}{3}$. This establishes the base case. Assume that $k\geq 1$
and each label tree $(T', S')\in \mathscr{T}$ with $h(T')<k$
satisfies the condition that $\frac{\gamma_t(T')}{\gamma_{sp}(T')}=
\frac{4}{3}$. Let $(T, S)\in \mathscr{T}$ be a label tree with
$h(T)=k$. Then $(T, S)$ can be obtained from a label tree $(T',
S')\in \mathscr{T}$ with $h(T')<k$ by the operation $\mathscr{O}$.
That is, add a path $u_1u_2u_3u_4u_5u_6$ and the edge $u_3v$, where
$v$ is a vertex of $T'$ which has status $B$. By induction,
$\frac{\gamma_t(T')}{\gamma_{sp}(T')}= \frac{4}{3}$.

By Observation~2.1(2), we can obtain a $\gamma_t$-set of $T$, say
$D$, which contains no leaf. Clearly, $\{v_2, v_3, v_4,
v_5\}\subseteq D$. Since $v$ has status $B$, it follows from
Observation~2.9(b) and (d) that $v$ has one neighbor which is a
support vertex of degree two, say $w$. Moreover, $\{v, w\}\subseteq
D$. It means that $D\setminus \{v_2, v_3, v_4, v_5\}$ is a total
dominating set of $T'$. That is, $\gamma_t(T)-4\geq \gamma_t(T')$.
On the other hand, let $D'$ be a $\gamma_t$-set of $T'$. It is easy
to see that $D'\cup \{v_2, v_3, v_4, v_5\}$ is a total dominating
set of $T$. That is, $\gamma_t(T')+4\geq \gamma_t(T)$. Therefore,
$\gamma_t(T')+4= \gamma_t(T)$.

Let $R'$ be a $\gamma_{sp}$-set of $T'$. By Proposition~2.4, there
exists a $\gamma_{sp}$-set $R$ of $T$ such that the leaf-neighbor of
$w$ belongs to $\overline{R}$. And then $w$ is its external private
neighbor with respect to $\overline{R}$. Moreover, $v\not\in
\overline{R}$. It implies that $R\cap V(T')$ is a super dominating
set of $T'$. Since $|R\cap \{u_1, u_2, u_3, u_4, u_5, u_6\}|\geq 3$,
we have that $|R|-3\geq |R'|$. On the other hand, let $K=R'\cup
\{v_2, v_3, v_6\}$ when $v\not \in \overline{R'}$ and $K=R'\cup
\{v_1, v_4, v_5\}$ when $v\in \overline{R'}$, then $K$ is a super
dominating set of $T$. That is, $|R'|+3\geq |R|$. Hence,
$\gamma_{sp}(T)=\gamma_{sp}(T')+3$.

In summary,
$\frac{\gamma_t(T)}{\gamma_{sp}(T)}=\frac{\gamma_t(T')+4}{\gamma_{sp}(T')+3}=\frac{4}{3}$.
\ep

Below we will prove Theorem~2.6.

\bp We proceed by induction on the order $n$ of the tree $T$(As
mentioned above, we only need to consider the case that $T$ has no
strong support vertex). The result is immediate for $n\leq 5$. Let
$n\geq 6$ and assume that for every tree $T'$ satisfying $3\leq
|T'|<n$, we have $\frac{\gamma_t(T')}{\gamma_{sp}(T')}\leq
\frac{4}{3}$, with equality if and only if $(T', S')\in \mathscr{T}$
for some labeling $S'$.

The result holds when $diam(T)\leq 5$. Moreover, if
$\frac{\gamma_t(T)}{\gamma_{sp}(T)}=\frac{4}{3}$, then $(T, S)=(P_6,
S_0)\in \mathscr{T}$. Hence, we may assume that $diam(T)\geq 6$. Let
$P=v_1v_2\cdots v_t$ be a longest path in $T$ such that $d(v_3)$ as
large as possible. We know that $d(v_2)=2$. Let $R$ be a
$\gamma_{sp}$-set of $T$ such that $v_1\in \overline{R}$. We now
proceed with two claims that we may assume are satisfied by the tree
$T$, for otherwise the desired result holds.

{\flushleft\textbf{Claim 1.}}\quad $d(v_3)=2$.

If not, assume that $d(v_3)>2$. Let $T_1=T-\{v_1, v_2\}$ and $D_1$
be a $\gamma_t$-set of $T_1$ which contains no leaf. Clearly,
$v_3\in D_1$, and $D_1\cup \{v_2\}$ is a total dominating set of
$T$. So, $\gamma_t(T)\leq \gamma_t(T_1)+1$. On the other hand, let
$R_1$ be a $\gamma_{sp}$-set of $T_1$. Note that $v_2$ is the
external private neighbor of $v_1$ with respect to $\overline{R}$.
And then $R\setminus \{v_2\}$ is a super dominating set of $T_1$.
So, $|R|-1\geq |R_1|$. By induction, we have that $3\gamma_t(T)\leq
3(\gamma_t(T_1)+1)=3\gamma_t(T_1)+3\leq 4\gamma_{sp}(T_1)+3\leq
4\gamma_{sp}(T)-1<4\gamma_{sp}(T)$. That is,
$\frac{\gamma_t(T)}{\gamma_{sp}(T)}<\frac{4}{3}$. \ep

{\flushleft\textbf{Claim 2.}}\quad $d(v_4)=2$.

Suppose that $d(v_4)>2$. Now we can distinguish three cases as
follows:

{\flushleft\textbf{Case 1.}}\quad $v_4$ is a support vertex.

Let $u$ be the leaf-neighbor of $v_4$, and $R'$ be a
$\gamma_{sp}$-set of $T$ such that $u\in \overline{R'}$. Then, $v_4$
is the external private neighbor of $u$ with respect to
$\overline{R'}$. It implies that $|\{v_1, v_2, v_3\}\cap
\overline{R'}|=1$ and $R'\setminus \{v_1, v_2, v_3\}$ is a super
dominating set of $T_1=T-\{v_1, v_2, v_3\}$. So, $|R'|-2\geq
|R'_1|$, where $R'_1$ is a $\gamma_{sp}$-set of $T_1$. On the other
hand, it is easy to see that $\gamma_t(T_1)+2\geq \gamma_t(T)$. By
induction, we have that $3\gamma_t(T)\leq
3(\gamma_t(T_1)+2)=3\gamma_t(T_1)+6\leq 4\gamma_{sp}(T_1)+6\leq
4\gamma_{sp}(T)-2<4\gamma_{sp}(T)$. That is,
$\frac{\gamma_t(T)}{\gamma_{sp}(T)}<\frac{4}{3}$.

{\flushleft\textbf{Case 2.}}\quad $v_4$ has a neighbor outside $P$,
say $u_1$, which is adjacent to a support vertex $u_2$.

From the choice of $P$, $d(u_1)=2$. Denote the leaf-neighbor of
$u_2$ by $u_3$. Note that $\{v_2, v_3\}\cap \overline{R}=\emptyset$.
If $v_4\in \overline{R}$, then $|\{u_1, u_2, u_3\}\cap
\overline{R}|=1$ and $R\setminus \{u_1, u_2, u_3\}$ is a super
dominating set of $T_1=T-\{u_1, u_2, u_3\}$. If $v_4\not \in
\overline{R}$, note that $|\{v_1, v_2, v_3\}\cap \overline{R}|=1$
and $R\setminus \{v_1, v_2, v_3\}$ is a super dominating set of
$T_1=T-\{v_1, v_2, v_3\}$. In either case, the proof is similar to
that of Case~1.

{\flushleft\textbf{Case 3.}}\quad Every neighbor of $v_4$ outside
$P$ is support vertex.

Let $u$ be a neighbor of $v_4$ outside $P$. Suppose that $d(v_4)\geq
4$, and let $T'$ be the component of $T-v_4u$ containing $v_4$.
Then, the proof is similar to that of Claim~1.

So we have that $d(v_4)=3$. Let $T''$ be the component of $T-v_4v_5$
containing $v_5$. It is easy to see that $\gamma_t(T'')+4\geq
\gamma_t(T)$. On the other hand, let $R_1$ be a $\gamma_{sp}$-set of
$T''$. Note that $\{v_2, v_3\}\cap \overline{R}=\emptyset$ and
$|\{u, w\}\cap \overline{R}|\leq 1$, where $w$ is the leaf-neighbor
of $u$. If neither $v_4$ nor $v_5$ belongs to $\overline{R}$, then
$(R\setminus \{v_2, v_3, w\})\cup \{v_1, u\}$ is a super dominating
set of $T$ whose cardinality is less than $R$, a contradiction. So,
$\{v_4, v_5\}\cap \overline{R}\neq \emptyset$.

If $v_4\in \overline{R}$, then $u\in \overline{R}$, and moreover,
$R\setminus \{v_2, v_3, w\})$ is a super dominating set of $T''$.
That is, $|R|-3\geq |R_1|$.

If $v_4\not \in \overline{R}$ and $v_5\in \overline{R}$, then
$|\{v_1, v_5, u, w\}\cap \overline{R}|=3$, and
$\overline{R}\setminus \{v_1, v_5, u, w\})$ is the complement of a
super dominating set of $T''$. That is, $|\overline{R}|-3\leq
|\overline{R_1}|$. Note that $|T|=|T''|+6$, we have that $|R|-3\geq
|R_1|$.

In either case, by induction, we have that $3\gamma_t(T)\leq
3(\gamma_t(T'')+4)=3\gamma_t(T'')+12\leq 4\gamma_{sp}(T'')+12\leq
4(\gamma_{sp}(T)-3)+12 =4\gamma_{sp}(T)$. That is,
$\frac{\gamma_t(T)}{\gamma_{sp}(T)}\leq \frac{4}{3}$. Suppose next
that $3\gamma_t(T)=4\gamma_{sp}(T)$. Then we have equality
throughout the above inequality chain. In particular,
$\gamma_t(T)=\gamma_t(T'')+4$ and $3\gamma_t(T'')=
4\gamma_{sp}(T'')$. By induction, $(T'', S')\in \mathscr{T}$ for
some labeling $S'$. If $v_5$ is a support vertex, by Lemma~2.10,
there exists a set $X$ of cardinality $\gamma_t(T'')-1$ such that
the leaf-neighbor of $v_5$, say $x$, belongs to $X$, and each vertex
of $T''$ is total dominated by $X$ except for $x$. In this case, let
$Y=(X\setminus \{x\})\cup \{v_5, v_2, v_3, v_4, u\}$. It is easy to
see that $Y$ is a total dominating set of $T$ with cardinality
$\gamma_t(T)-1$, it is impossible. If $v_5$ is a leaf, we can obtain
a contradiction through the similar argument. Hence, $v_5$ is
neither a leaf nor a support vertex. It follows from
Observation~2.9(a) and (b) that $v_5$ has status $B$ in $(T'', S')$.
Let $S$ be obtained from the labeling $S'$ by labeling the vertices
$v_1$ and $w$ with label $C$, the vertices $v_2$ and $u$ with label
$A$, the vertex $v_3$ and $v_4$ with label $B$. Then, $(T, S)$ can
be obtained from $(T'', S')$ by operation $\mathscr{O}$. Thus, $(T,
S)\in \mathscr{T}$.\ep

In summary, $d(v_2)=d(v_3)=d(v_4)=2$. Let $T'$ be the component of
$T-v_4v_5$ containing $v_5$, and $R_1$ be a $\gamma_{sp}$-set of
$T'$. It is easy to see that $\gamma_t(T')+2\geq \gamma_t(T)$. On
the other hand, similar to Case~3 of Claim~2, we have that
$|R|-2\geq |R_1|$. It follows that $3\gamma_t(T)\leq
3(\gamma_t(T')+2)=3\gamma_t(T')+6\leq 4\gamma_{sp}(T')+6\leq
4\gamma_{sp}(T)-8+6=4\gamma_{sp}(T)-2<4\gamma_{sp}(T)$. That is,
$\frac{\gamma_t(T)}{\gamma_{sp}(T)}<\frac{4}{3}$. \ep

\section{Bound on the super domination subdivision number of trees}

In this section, we first present the upper bound of
$sd_{\gamma_{sp}}(T)$.

\begin{thm}
For any tree $T$ of order at least $2$, $sd_{\gamma_{sp}}(T)\leq 2$.
\end{thm}

\bp It is easy to see that the result holds for a tree of
$diam(T)\leq 3$, so we assume that $diam(T)\geq 4$. Let
$P=u_1u_2u_3\cdots u_t$ be a longest path of $T$. Let $T'$ be
obtained from $T$ by subdividing the edges $u_2u_3$ and $u_3u_4$
with vertices $x$ and $y$. By Proposition~2.4, there exists a
$\gamma_{sp}$-set of $T'$, say $S'$, such that $\overline{S'}$
contains the vertex $u_1$. Observe that $u_2, x\not \in
\overline{S'}$. Let $D=(\overline{S'}\setminus \{u_1, y\})\cup
\{u_2\}$ when $u_3, y\in \overline{S'}$, $D=\overline{S'}\setminus
\{u_3, y\}$ when $|\{u_3, y\}\cap \overline{S'}|=1$,
$D=\overline{S'}\setminus \{u_4\}$ when $u_3, y\not \in
\overline{S'}$ and $u_4\in \overline{S'}$, $D=\overline{S'}$ when
$u_3, y, u_4\not \in \overline{S'}$. In either case, we note that
$V(T)\setminus D$ is a super dominating set of $T$. Combining the
fact that $|T'|=|T|+2$, we have that $\gamma_{sp}(T')\geq
\gamma_{sp}(T)+1$. That is, $sd_{\gamma_{sp}}(T)\leq 2$. \ep

Trees are classified as Class~1 or Class~2 depending on whether
their super domination subdivision number is 1 or 2, respectively.
Next, we are ready to provide the constructive characterizations of
trees in Class~1 and Class~2. We introduce the operation as follows.

{\bf Operation $\mathscr{U}_1$}: Add a star of order at least two
and join its center vertex to a vertex $v$ of $T'$ when there exists
a $\gamma_{sp}$-set $S$ of $T'$ such that $N_{T'}[v]\cap
\overline{S}=\emptyset$, or $v\not \in \overline{S}$ and
$N_{T'}[v]\cap U_S(T')=\emptyset$.

We define the family $\mathscr{U}$ as:

$\mathscr{U}=\{T|T$ is obtained from a star of order at least three
by a finite sequence of operation $\mathscr{U}_1\}$. We show first
that every tree in the family $\mathscr{U}$ is in Class~2.

\begin{lem}
If $T\in \mathscr{U}$, then $T$ is in Class~2.
\end{lem}

\bp The proof is by induction on the number $h(T)$ of operations
required to construct the tree $T$. Observe that $T$ is a star of
order at least three when $h(T)=0$, and the result holds. This
establishes the base case. Assume that $k\geq 1$ and each tree
$T'\in \mathscr{U}$ with $h(T')<k$ is in Class~2. Let $T\in
\mathscr{U}$ be a tree with $h(T)=k$. Then $T$ can be obtained from
a tree $T'\in \mathscr{U}$ with $h(T')<k$ by the operation
$\mathscr{U}_1$. In other words, $T$ is obtained from $T'$ by adding
a star of order at least two and join its center vertex, say $u$, to
a vertex $v$ of $T'$, where $N_{T'}[v]\cap \overline{S'}=\emptyset$,
or $v\not \in \overline{S'}$ and $N_{T'}[v]\cap
U_{S'}(T')=\emptyset$, $S'$ is some $\gamma_{sp}$-set of $T'$. By
induction, $T'$ is in Class~2.

Let $u_1$ be a leaf-neighbor of $u$ in $T$, and $S$ be a
$\gamma_{sp}$-set of $T$ such that $u_1\in \overline{S}$. Then, $u$
is the external private neighbor of $u_1$ with respect to
$\overline{S}$, and it means that $\overline{S}\setminus \{u_1\}$ is
the complement of a super dominating set of $T'$, and so
$|\overline{S}|-1\leq |\overline{S'}|$. Moreover, note that $v\not
\in \overline{S'}$, $\overline{S'}\cup \{u_1\}$ is the complement of
a super dominating set of $T$, so $|\overline{S'}|+1\leq
|\overline{S}|$. Hence, $|\overline{S'}|+1=|\overline{S}|$.

Let $e\in E(T)$ and $T^{*}$ be obtained from $T$ by subdividing the
edge $e$ with vertex $x$. Let $S^{*}$ be a $\gamma_{sp}$-set of
$T^{*}$. Now we can distinguish three cases as follows:

{\flushleft\textbf{Case 1.}}\quad $e\in E(T')$.

Let $T'^{*}$ be obtained from $T'$ by subdividing the edge $e$ with
a vertex, and $S'^{*}$ be a $\gamma_{sp}$-set of $T'^{*}$. Similar
to the argument as above, we have that
$|\overline{S'^{*}}|+1=|\overline{S^{*}}|$. On the other hand, by
induction, $\gamma_{sp}(T'^{*})=\gamma_{sp}(T')$. And then,
$|\overline{S'}|+1=|\overline{S'^{*}}|$. It concludes that
$|\overline{S'}|+2=|\overline{S^{*}}|$. \ep

{\flushleft\textbf{Case 2.}}\quad $u$ is a strong support vertex in
$T$ and $e=uu_1$.

Suppose that $u_2$ is a leaf-neighbor of $u$ in $T$ other than
$u_1$. Note that $v\not \in \overline{S'}$, then $\overline{S'}\cup
\{u_1, u_2\}$ is the complement of a super dominating set of
$T^{*}$. That is, $|\overline{S'}|+2\leq |\overline{S^{*}}|$. On the
other hand, by Proposition~2.4, without loss of generality, we may
assume that $u_2\in \overline{S^{*}}$. Then, $u$ is the external
private neighbor of $u_2$ with respect to $\overline{S^{*}}$.
Moreover, we have that $u_1\in \overline{S^{*}}$ and $x$ is the
external private neighbor of $u_1$ with respect to
$\overline{S^{*}}$. Therefore, $\overline{S^{*}}\setminus \{u_1,
u_2\}$ is the complement of a super dominating set of $T'$. That is,
$|\overline{S^{*}}|-2\leq |\overline{S'}|$. So, we have that
$|\overline{S^{*}}|-2=|\overline{S'}|$. \ep

{\flushleft\textbf{Case 3.}}\quad $u$ is a strong support vertex in
$T$ and $e=uv$, or $u$ is not a strong support vertex in $T$ and
$e\in \{uv, uu_1\}$.

We assume that $u$ is not a strong support vertex in $T$ and
$e=uu_1$(The other two cases can also be discussed similarly). Let
$D=\overline{S'}\cup \{u_1, v\}$ when $v\not \in \overline{S'}$ and
$N_{T'}[v]\cap U_{S'}(T')=\emptyset$, and $D=\overline{S'}\cup \{u,
x\}$ when $N_{T'}[v]\cap \overline{S'}=\emptyset$. $D$ is the
complement of a super dominating set of $T^{*}$. That is,
$|\overline{S'}|+2\leq |\overline{S^{*}}|$. On the other hand, by
Proposition~2.4, without loss of generality, we assume that $u_1$ is
in $\overline{S^{*}}$. And then $x$ is the external private neighbor
of $u_1$ with respect to $\overline{S^{*}}$. Among all vertices of
$\overline{S^{*}}\setminus \{u_1\}$, let $y$ be the vertex at
minimum distance from $u$. It is easy to see that
$\overline{S^{*}}\setminus \{u_1, y\}$ is the complement of a super
dominating set of $T'$. That is, $|\overline{S^{*}}|-2\leq
|\overline{S'}|$. Hence, $|\overline{S^{*}}|-2=|\overline{S'}|$. \ep

In either case, we have that $|\overline{S^{*}}|-2=|\overline{S'}|$.
Combining the fact that $|\overline{S'}|+1=|\overline{S}|$, we have
that $|\overline{S^{*}}|=|\overline{S}|+1$. It follows from
$|T^{*}|=|T|+1$ that $\gamma_{sp}(T^{*})=\gamma_{sp}(T)$. That is,
$T$ is in Class~2.
 \ep

\begin{lem}
If a tree $T$ is in Class~2, then $T\in \mathscr{U}$.
\end{lem}

\bp We know that $T$ is in Class~2, it is a star of order at least
three when $diam(T)\leq 2$, and $T\in \mathscr{U}$. So we consider
the case when $diam(T)\geq 3$. We proceed by induction on the order
$n$ of $T$. Assume that the result is true for every tree in Class~2
of order less than $n$. Let $P=v_1v_2\cdots v_t$ be a longest path
in $T$, and $T_1$ be the component of $T-v_2v_3$ containing $v_3$.
Let $e\in E(T_1)$ and $T_1^{*}$ (respectively, $T^{*}$) be obtained
from $T_1$ (respectively, $T$) by subdividing the edge $e$, and $S$
(respectively, $S^{*}$, $S_1$, $S_1^{*}$) be a $\gamma_{sp}$-set of
$T$ (respectively, $T^{*}$, $T_1$, $T_1^{*}$).

Set $D=\overline{S_1}\cup \{v_1\}$ when $v_3\not \in \overline{S_1}$
and $D=\overline{S_1}\cup \{v_2\}$ when $v_3\in \overline{S_1}$. It
is easy to see that $D$ is the complement of a super dominating set
of $T$. That is, $|\overline{S_1}|+1\leq |\overline{S}|$. On the
other hand, without loss of generality, assume that $v_1\in
\overline{S}$, then $\overline{S}\setminus \{v_1\}$ is the
complement of a super dominating set of $T_1$. And so,
$|\overline{S}|-1\leq |\overline{S_1}|$. Hence,
$|\overline{S}|-1=|\overline{S_1}|$. Similarly, we have that
$|\overline{S^{*}}|-1=|\overline{S_1^{*}}|$.

By assumption, we know that $\gamma_{sp}(T)=\gamma_{sp}(T^{*})$.
That is, $|\overline{S^{*}}|-1=|\overline{S}|$. So,
$|\overline{S_1}|=|\overline{S}|-1=|\overline{S^{*}}|-2=|\overline{S_1^{*}}|-1$.
It implies that $T_1$ is in Class~2. By induction, $T_1\in
\mathscr{U}$.

Now, let $T'$ be the tree obtained from $T$ by subdividing the edge
$v_2v_3$ with vertex $x$, $T''$ be the component of $T'-xv_3$
containing $x$, and let $S'$ be a $\gamma_{sp}(T')$-set. According
to the discussion as above, we know that
$|\overline{S'}|=|\overline{S}|+1$ and
$|\overline{S_1}|=|\overline{S}|-1$. It concludes that
$|\overline{S'}|=|\overline{S_1}|+2$. If $|V(T'')\cap
\overline{S'}|\geq 2$, then we have that $\{x, v_2\}\subseteq
\overline{S'}$. It means that $v_3$ is the external private neighbor
of $x$ with respect to $\overline{S'}$. Moreover, $(N(v_3)\setminus
\{x\})\cap \overline{S'}=\emptyset$. Let $D=\overline{S'}\setminus
\{x, v_2\}$, $V(T_1)\setminus D$ is a super dominating set of $T_1$.
Since $|\overline{S'}|=|\overline{S_1}|+2$, $V(T_1)\setminus D$ is a
$\gamma_{sp}(T_1)$-set. And $T$ can be obtained from $T_1$ by
operation $\mathscr{U}_1$.

If $|V(T'')\cap \overline{S'}|=1$, without loss of generality,
assume that $v_1\in \overline{S'}$. In this case, $v_2, x\not \in
\overline{S'}$, and $(N[v_3]\setminus \{x\})\cap \overline{S'}\neq
\emptyset$. Now we can distinguish three cases as follows:

{\flushleft\textbf{Case 1.}}\quad $v_3\in \overline{S'}$.

Let $D=\overline{S'}\setminus \{v_1, v_3\}$. Note that
$|\overline{S'}|=|\overline{S_1}|+2$, $H=V(T_1)\setminus D$ is a
$\gamma_{sp}(T_1)$-set satisfying $v_3\not \in D$ and
$N_{T_1}[v_3]\cap U_H(T_1)=\emptyset$. That is, $T$ can be obtained
from $T_1$ by operation $\mathscr{U}_1$.

{\flushleft\textbf{Case 2.}}\quad $v_3\not \in \overline{S'}$ and
$v_3\in P_{S'}(T')$.

Since $v_3\in P_{S'}(T')$, there exists a $y\in (N(v_3)\setminus
\{x\})\cap \overline{S'}$ such that $v_3$ is the external private
neighbor of $y$ with respect to $\overline{S'}$. Let
$D=(\overline{S'}\setminus \{y, v_1\})\cup \{x, v_2\}$. Clearly,
$V(T')\setminus D$ is also a $\gamma_{sp}(T')$-set and the proof is
similar to the case of $|V(T'')\cap \overline{S'}|\geq 2$.

{\flushleft\textbf{Case 3.}}\quad $v_3\not \in \overline{S'}$ and
$v_3\not \in P_{S'}(T')$.

In this case, $|(N(v_3)\setminus \{x\})\cap \overline{S'}|\geq 2$.

{\flushleft\textbf{subcase 3.1.}}\quad $(N(v_3)\setminus \{x\})\cap
P_{S'}(T')=\emptyset$.

Take a vertex $z\in (N(v_3)\setminus \{x\})\cap \overline{S'}$. Let
$D=(\overline{S'}\setminus \{z\})\cup \{v_3\}$. Note that the proof
is similar to the case~1.

{\flushleft\textbf{subcase 3.2.}}\quad $(N(v_3)\setminus \{x\})\cap
P_{S'}(T')\neq \emptyset$.

Assume that $(N(v_3)\setminus \{x\})\cap P_{S'}(T')=\{u_1, u_2,
\cdots, u_k\}$. For any $i\in \{1, 2, \cdots, k\}$, let $T'_i$ be
the component of $T'-v_3u_i$ containing $u_i$. Clearly, $S'_i=S'\cap
V(T'_i)$ is a $\gamma_{sp}(T'_i)$-set. By Proposition~2.4, for each
$i\in \{1, 2, \cdots, k\}$, there must be a $\gamma_{sp}(T'_i)$-set
$S''_i$ such that $u_i\in \overline{S''_i}$. Then, $(S'\setminus
(\bigcup \limits_{i=1}^{k} S'_i))\cup (\bigcup \limits_{i=1}^{k}
S''_i)$ is a $\gamma_{sp}(T')$-set and the proof is similar to the
subcase 3.1. \ep

As an immediate consequence of Lemmas~3.2 and 3.3 we have the
following theorem.

\begin{thm}
A tree $T$ is in Class~2 if and only if $T\in \mathscr{U}$.
\end{thm}

Let $\mathscr{G}=\{T|T$ is a nontrivial tree$\}$, and
$\mathscr{P}=\mathscr{G}\setminus \mathscr{U}$. We immediately
obtain the following result.

\begin{thm}
A tree $T$ is in Class~1 if and only if $T\in \mathscr{P}$.
\end{thm}


\end{document}